\documentclass[12pt]{amsart}
\usepackage{amscd,amssymb}
\usepackage[arrow,matrix]{xy}
\usepackage[plainpages,backref,urlcolor=blue]{hyperref}

\theoremstyle{plain}

\newtheorem{prop}[subsection]{Proposition}
\newtheorem{cor}[subsection]{Corollary}

\theoremstyle{definition}
\newtheorem{rk}[subsection]{Remark}
\newtheorem{definition}[subsection]{Definition}
\newtheorem{ex}[subsection]{Example}

\numberwithin{equation}{section}
\newcommand{\NN}{{\mathcal N}}

\newcommand{\C}{\mathbb{C}}

\newcommand{\PP}{\mathbb{P}}

\DeclareMathOperator{\defect}{def}
\DeclareMathOperator{\codim}{codim}

\begin{document} 

\title [Invariants and rigidity of projective hypersurfaces]
{Invariants and rigidity of projective hypersurfaces} 

\author[Gabriel Sticlaru]{Gabriel Sticlaru}
\address{Faculty of Mathematics and Informatics,
Ovidius University,
Bd. Mamaia 124, 900527 Constanta,
Romania}
\email{gabrielsticlaru@yahoo.com }



\begin{abstract}

 This paper continues our researches \cite{DS1, DS2, DS3} by computing some invariants based on Hilbert-Poincar\'{e} series associated to Milnor algebras. Our computations are for some of the classical surfaces  and 3-folds with different configurations of isolated singularities.
As a by-product of a recent result of E. Sernesi, we give examples of classical hypersurfaces which are (or are not) projectively rigid.
We also include a Singular program to compute the invariants and to decide if a singular projective hypersurface is nodal and projectively rigid.  
 
\end{abstract}

\maketitle
 
\section{Introduction}

Let $S=\C[x_0,...,x_n]$ be the graded ring of polynomials in $x_0,,...,x_n$ with complex coefficients and denote by $S_r$ the vector space of homogeneous polynomials in $S$ of degree $r$. 
For any polynomial $f \in S_r$ we define the {\it Jacobian ideal} $J_f \subset S$ as the ideal spanned by the partial derivatives $f_0,...,f_n$ of $f$ with respect to $x_0,...,x_n$. For $n=2$ we use $x,y,z$ instead of
$x_0, x_1, x_2$ and $f_x, f_y, f_z$  instead of $f_0, f_1, f_2$.

The Hilbert-Poincar\'{e} series of a graded $S$-module $M$ of finite type is defined by 
\begin{equation} 
\label{eq2}
HP(M)(t)= \sum_{k\geq 0} \dim M_k\cdot t^k 
\end{equation} 
and it is known, to be a rational function of the form 
\begin{equation} 
\label{eq3}
HP(M)(t)=\frac{P(M)(t)}{(1-t)^{n+1}}=\frac{Q(M)(t)}{(1-t)^{d}}.
\end{equation}

For any polynomial $f \in S_r$ we define the corresponding graded {\it Milnor} (or {\it Jacobian}) {\it algebra} by
\begin{equation} 
\label{eq1}
M=M(f)=S/J_f.
\end{equation}
In fact, such a Milnor algebra can be seen (up to a twist in grading) as the first (or the last) homology (or cohomology) of the Koszul complex of the partial derivatives $f_0,...,f_n$ in $S$, see  \cite{CD} or  \cite{D1}, Chapter 6.

One of our research aims will be to improve the bounds in Choudary-Dimca Theorem from \cite {CD}, to get sharp estimates in many cases.

\textbf{Choudary-Dimca Theorem }
Let $V(f): f=0$ a hypersurface in $ \PP^n$ with only isolated singularities. 
For any $q \geq T+1, T=(n+1)(d-2)$, one has 
$$ \dim M(f)_q= \tau(V(f))=\sum_{j=1,p}\tau(V(f),a_j) $$
where $\tau(V(f))$ is the global Tjurina number of the hypersurface $V(f)$. In particular, the
Hilbert polynomial $H(M(f))$ is constant and this constant is $\tau(V(f))$.

\bigskip
For a hypersurface $D: f=0$ in $\PP^n$ with isolated singularities we recall \textbf{four} invariants, introduced in 
\cite{DS2} and add a new one, as follows.

\begin{definition}
\label{def}

\noindent (i) The {\it coincidence threshold} $ct(D)$ defined as
$$ct(D)=\max \{q~~:~~\dim M(f)_k=\dim M(f_s)_k \text{ for all } k \leq q\},$$
with $f_s$  a homogeneous polynomial in $S$ of degree $d=\deg f$ such that $D_s:f_s=0$ is a smooth hypersurface in 
$ \PP^n$.

\bigskip

\noindent (ii) The {\it stability threshold} $st(D)$ defined as
$$st(D)=\min \{q~~:~~\dim M(f)_k=\tau(D) \text{ for all } k \geq q\}$$
where $\tau(D)$ is the total Tjurina number of $D$, i.e. the sum of all the Tjurina numbers of the singularities of $D$.

\bigskip

\noindent (iii) The {\it minimal degree of a nontrivial syzygy} $mdr(D)$ defined as
$$mdr(D)=\min \{q~~:~~ H^n(K^*(f))_{q+n}\ne 0\}$$
where $K^*(f)$ is the Koszul complex of $f_0,...,f_n$ with the natural grading.

\bigskip

\noindent (iv)  
Let $D_{smooth}:f_{smooth}=0$ be a smooth hypersurface of the same degree $d$ in $ \PP^n$. 
We define the integer $def(D)=$defect of $D$ as 
$$ \text {def(D) = the first not zero coefficient of the difference } S(t)-F(t)$$
where $S(t)$ (resp. $F(t)$ are the corresponding  Hilbert-Poincar\'{e} series of Milnor Algebras
$M(f)$ (resp. $M(f_{smooth}$).

\bigskip

\noindent (v) 
A finite sequence of strictly positive real numbers $a_0,...,a_q$ is said to be log-concave if
$a_{i-1}a_{i+1}\leq a_i^2$ for all $i=1,2,...,q-1$.

We define a new integer $lc(D)=$  {\it the  log-concavity}  of $D$ as the maximal integer $q$ such that the sequence 
$ \dim M(f)_0, \ldots , \dim M(f)_q $ is  a log-concave sequence. If $q \geq st(D)+1$ we put $lc(D)=\infty. $ 

\end{definition} 
Recall also that, for a finite set of points $\NN \subset \PP^n$,
we denote by 
$$\defect S_m(\NN)=|\NN| - \codim \{h \in S_m~~|~~ h(a)=0 \text{ for any } a \in \NN\},$$
the {\it defect (or superabundance) of the linear system of polynomials in $S_m$ vanishing at the points in $\NN$}, see \cite{D1}, p. 207. This positive integer is called the {\it failure of $\NN$ to impose independent conditions on homogeneous polynomials of degree $m$}. \\

Except for the explicit computations and the final Singular program, our main results are 
Proposition \ref{prop2} and Corollary \ref{cor1} and \ref{cor2} relating our invariants to the question of projective rigidity of hypersurfaces with isolated singularities.

\section {Computations for some projective hypersurfaces}

In this section, in the first part we analyze  two curves, of degree 4, one nodal and one cuspidal, in order to explain in detail our approach. \\
For each curve $C: 	f=0 $, we find and classify all the singularities, give the genus and compute 
the Hilbert-Poincar\'{e} series and our invariants.

Firstly, we will find all multiple points of $C : f(x,y,z)=0$ by solving the system of equations: \\
$f(x,y,z)=0, \frac{\partial f}{\partial x}(x,y,z) =0, \frac{\partial f}{\partial y}(x,y,z) =0, \frac{\partial f}{\partial z}(x,y,z) =0 $ \\
Because f is homogenous of degree $d$, we have:\\
$x \frac{\partial f}{\partial x}(x,y,z)+ y \frac{\partial f}{\partial y}(x,y,z) + 
z \frac{\partial f}{\partial z}(x,y,z) = d f(x,y,z).$ \\
Obviously, if all partial derivatives of polynomial f vanish at $p$, then the
polynomial $f$ vanishes at $p$ too. Therefore, for finding singular points it is enough to solve the system of equations:\\
$  \frac{\partial f}{\partial x}(x,y,z) =0, \frac{\partial f}{\partial y}(x,y,z) =0, \frac{\partial f}{\partial z}(x,y,z) =0. $ \\

In this computation we will use the $\delta$-invariant of a plane curve singularity. To determine this invariant, the simplest way is to use the Milnor-Jung formula:
\begin{equation} 
\label{MJ}
\mu=2 \delta-r+1,
\end{equation}
where $\mu$ is the Milnor number and $r$ the number of branches of the singularity. For example, for the $A_3$ singularity $x^2-y^4=0$ we have $\mu=3$, $r=2$ and hence $\delta= 2$.

The Hilbert-Poincar\'{e} series is computed with two methods: combinatorial and based on a free resolution.

To compute invariants $ct $ and $def ,$ 
we need the Hilbert-Poincar\'{e} series $F(t)$ for $d$ degree Fermat curve:  $x^d+y^d+z^d=0$ 
(the same for any smooth $d$ degree curve):\\ 
$F_d(t)=(1+t+\ldots+t^{d-2})^3$ and for $d=4,$ $F_4(t)=1+3t+6t^2+7t^3+6t^4+3t^5+t^6. $\\

It is a very hard work to obtain manually these informations, even if we consider curves and surfaces with low degree. 

In the secont part, all the compuations was made by our Singular program.

\subsection {Lemniscate of Bernoulli and Cardioid}

\noindent 
The Lemniscate of Bernoulli and the Cardioid are splendid curves of the fourth degree, with different topological type of singularity (but with the same singular points).
 
The \textbf{Lemniscate of Bernoulli} is a nodal curve with $3$ nodes (type $3A_1$). 
The affine equation of the Lemniscate of Bernoulli looks like 
$$F(x,y)=(x^2+y^2)^2- 2(x^2-y^2).$$  

Next, we will homogenize the defining polynomial  and we get the defining polynomial
of its associated projective curve:
$$ L: f(x,y,z)=(x^2+y^2)^2-2(x^2-y^2)z^2=0 .$$ 

\textit{Find the position of the singularities, i.e. solve the system:} \\
$ \frac{\partial f}{\partial x}(x,y,z) = 4x(x^2+y^2)-4xz^2=0 $\\ 
$ \frac{\partial f}{\partial y}(x,y,z) = 4yz^2+4y(x^2+y^2)=0 $ \\ 
$ \frac{\partial f}{\partial z}(x,y,z) = -4(x^2-y^2)z=0 $ \\
For $z = 1$ we obtain the double point a=(0:0:1).
For $z = 0$ (singularities on the line at infinity $L_\infty  : z=0)$ we get two complex double points $b=(1:i:0)$ and $c=(1:-i:0).$ Hence $3$ singularities $a,b,c$ with multiplicity $2.$ 

\textit{Determine the type of singularities }
We start with $a=(0:0:1).$\\
With local coordinates at a: $u=\frac{x}{z},v=\frac{y}{z}$, local equation at a is:\\ 
$g=2(u^2-v^2)-(u^2+v^2)^2=0$. \\
Use the weights $w_1=wt(u)=1$ and $w_2=wt(v)=1$ \\
Then $g$ is semi-weighted homogenous of type $(w_1,w_2,d)$ with $d=2$ \\
Prop (7.37), p 116 and Prop (7.27), p 112 \cite{D3} imply that $\mu (g)= \mu (u^2-v^2)=1$ so 
$ g \equiv A_1 $ is a node (the only singularity with $\mu=1).$ \\
For nodes it is known that 
$ \delta (A_1) =1 $ and $r(A_1)=2 $ (two branches $ u-v=0 $ and $ u+v=0).$ \\
\emph{Let's treat now $ b=(1:i:0) $} \\
Local coordinates $ u= \frac{x+iy}{x}, v=\frac{z}{x}. $\\
$ \frac{y}{x}= -i(u-1) \Leftrightarrow u=1+ i\frac{y}{x}, v=\frac{z}{x}. $\\
we divide by $x^4 :$ $(1+ (\frac{y}{x})^2 )^2 - 2(1-( \frac{y}{x})^2 )(\frac{z}{x})^2 = 0 $ \\
Local equation at b:
$ g=\left[ 1-(u-1)^2 \right] ^2 - 2(1+(u-1)^2)v^2 = \\
= (1-u^2+2u-1)^2-2v^2-2v^2(1-2u+u^2)=\\
=4u^2-4v^2+ $ higher degree terms.\\
Hence with the same weight as before g is semi-weighted homogenous and \\
$ \mu (g) = \mu (4(u^2-v^2))=1 \Rightarrow g \equiv A_1 .$\\
Same works for c.  
\textit{ Computation of genus}\\
Genus of a smouth curve of degree $d=4$ is:\\
$ g_s =\frac{(d-1)(d-2)}{2}=\frac{6}{2}=3. $ \\
Genus of our singular curve is: \\
$ g(L) = g_s - \sum _{x \in \left\{ a,b,c \right\} } \delta (x) = 3-3 = 0. $ Hence $L$ is rational (we have to check also that f is irreducible i.e. $ f \neq f_1 (x,y,z) \times f_2(x,y,z)$ with $f_1 , f_2$ homogenous of degree $ d_1 > 0, d_2 > 0, d_1+d_2 = 4). $ \\
Hence, the genus of $L$ is $0$ and $L$ is irreducible. Thus, the curve $L$ is rational parameterizable and:\\
$x(t)=2(-3 - 2t + 2t^3 + 3t^4)/(5 + 12t + 30t^2 + 12t^3 + 5t^4)$ \\
$y(t)=-2(-1 - 6t + 6t^3 + t^4)/(5 + 12t + 30t^2 + 12t^3 + 5t^4) $\\
is a rational parametrization for the affine equation.\\
To compute the Hilbert-Poincar\'{e} series, first recall that  the quotient rings $S/I$ and $S/LT(I)$ have the same series, for any monomial ordering, where LT(I) is the ideal of leading terms of the ideal I. The leading ideal of the jacobian $J_f $ is:

$ LI=<yz^4, y^2z^3, y^3z, xz^4, xyz^2, xy^2z, x^2z, x^2y, x^3>$

For the graded Milnor algebra, $M=\oplus_{k\geq0}M_{k} $ we show the bases for the  homogeneous components: $ M_0= \left\{  1 \right\}, M_1= \left\{   z,y,x  \right\},\\
M_2= \left\{   z^2,yz,xz,y^2,xy,x^2  \right\},M_3= \left(  z^3,yz^2,xz^2,y^2z,xyz,y^3,xy^2  \right\},\\
M_4= \left(  z^4,yz^3,xz^3,y^2z^2,y^4,xy^3  \right\}$ and  $M_k= \left\{   z^k,y^k,xy^{k-1}  \right\}$ for all $k\geq 5$.

Finally, if we count the number of monomials in each homogeneous bases, we find the Hilbert-Poincar\'{e} series 
$ S(t)=1+3t+6t^2+7t^3+6t^4+3(t^5+\ldots $ \\
\textit{Our invariants are:} $\tau=3,\  \  ct=5, \  \  st=5, \  \  mdr=3, \  \  def=2 \  \ lc=5 $
and because $mdr=3$, we show three nontrivial  linear independent relations between derivatives (syzygies), with polynomial coefficients  of degree 3:\\
$
(2y^3+yz^2)f_x + (xz^2 - 2xy^2)f_y +2xyzf_z=0 \\
(2xy^2+xz^2)f_x + (yz^2-2x^2y)f_y +(x^2z -y^2z-z^3)f_z =0\\
(2x^2y+yz^2)f_x+(xz^2-2x^3)f_y-2xyzf_z =0 \\
$ 

The \textbf{Cardioid} is a cuspidal curve, with $3$ cusps, (type $3A_2$) \\
The affine equation of Cardioid is:\\
$ F(x,y)=(x^2+y^2+x)^2 - (x^2+y^2) =0 $ 

The projective ecuation is:
$ C: f(x,y,z)=(x^2+y^2+xz)^2 - (x^2+y^2)z^2=0. $ 

\textit{Find the position of the singularities, i.e. solve the system:} \\
$ \frac{\partial f}{\partial x}(x,y,z) = 4x(x^2+y^2+xz)+2x^2z+2y^2z=0 $\\ 
$ \frac{\partial f}{\partial y}(x,y,z) = 2y(2x^2+2y^2+2xz-z^2)=0 $ \\ 
$ \frac{\partial f}{\partial z}(x,y,z) = 2(x^2+y^2+xz)x-2z(x^2+y^2)=0 $ \\
For $y = 0$ we obtain the double point $a=(0:0:1).$
For $z = 0$ we get two complex double points $b=(1:i:0)$ and $c=(1:-i:0).$ Hence 3 singularities $a,b,c.$ 

\textit{Determine the type of singularities }\\
Consider the case $a=(0:0:1).$\\
With local coordinates at a: $u=\frac{x}{z},v=\frac{y}{z}$, local equation at a is:\\ 
$g=2(u^2+v^2+u)^2-(u^2+v^2)= 2u^3-v^2+(u^4+v^4+2u^2v^2+2uv^2).$. \\
Use the weights $w_1=wt(u)=2$ and $w_2=wt(v)=3$ \\
Then $g$ is semi-weighted homogenous of type $(w_1,w_2,d)$ with $d=6$ \\
Prop (7.37), p 116 and Prop (7.27), p 112 \cite{D3} imply that $\mu (g)= \mu (2u^3-v^2)=2$ so 
$ g \equiv A_2 $ is a cusp \\
Similar computations for $b$ and $c$. It is known that $ \delta (A_2) =1 $ and $r(A_2)=1 $ \\
Computation genus as before, 
$ g(C) = g_s - \sum _{x \in \left\{ a,b,c \right\} } \delta (x) = 3-3 = 0. $ 
Hence, the genus is 0, since the Cardioid is irreducible. Thus, the curve $C$ is rational parameterizable and:\\
$
x(t)=2(-1 + 4t^2)/(1 + 4t^2)^2 \\
y(t)=-8t/(1 + 4t^2)^2
$
is a rational parametrization for the affine equation.\\
Here is the minimal graded free resolution of Milnor algebra M:
\begin{equation} 
\label{res}
0 \to R_3 \stackrel{C}{\rightarrow} R_2 \stackrel{B}{\rightarrow} R_1 \stackrel{A}{\rightarrow} R_0 \to M \to 0
\end{equation}
where $R_0=S$, $R_1=S^3(-3)$, $R_2=S^{3}(-5)$ and $R_3= S(-6)$ and the matrix are: 

{\em A is}  $1 \times 3$~matrix
    \[ \left( \begin{array}{ccc}
        2x^3 + 2xy^2 - 2y^2z   & 4x^2y + 4y^3 + 4xyz - 2yz^2  & 6x^2z + 6y^2z  
     \end{array} \right)\] 
   {\em B is}
    $3 \times 3$~matrix
    \[ \left( \begin{array}{ccc}
    -2y^2 - 1/2z^2  & -2xy + yz & 1/2xz - 1/4z^2 \\
    xy + 1/2yz & x^2 & 1/4yz \\
    -y^2 + 1/6xz & -xy + 1/3yz & -1/6x^2 - 1/6y^2 + 1/12xz \end{array} \right)\] 
   {\em C is}
    $3 \times 1$~matrix
    \[ \left( \begin{array}{ccc}
      -1/2x + 1/4z     \\
          1/2y       \\
          -1/2z         \end{array} \right)\] 
To get the formulas for the Hilbert-Poincar\'{e} series, we start with the resolution \eqref{res} and get
$HP(M)(t)=HP(R_0)(t)-HP(R_1)(t)+HP(R_2)(t)-HP(R_3)(t).$

Then we use the well-known formulas $HP(N\oplus N')(t)=HP(N)(t)+HP( N')(t)$, $HP(N(-r))(t)=t^rHP(N)(t)$
and $HP(S)(t)=\frac{1}{(1-t)^3}.$ and we obtain:

$$HP(M)(t)=\frac{1-3t^3+3t^5-t^6}{(1-t)^3}=\frac{1+2t+3t^2+t^3-t^4}{1-t}=1+3t+6t^2+7t^3+6(t^4+\ldots $$

\textit{Our invariants are:}
$\tau=6, \  \  ct=4, \  \  st=4, \  \  mdr=2, \  \  def=3 \  \ lc=4. $ and because $mdr=2$, we show three nontrivial  linear independent relations between derivatives (syzygies), with polynomialcoefficients of degree 2:\\
$
(xz-6y^2)f_x +(6xy+3yz)f_y-(3z^2+2xz)f_z =0\\
(yz-3xy)f_x + 3x^2f_y+yzf_z =0\\
(xz-2x^2-2y^2)f_x +3yzf_y +(4x^2+4y^2+4xz-3z^2)f_z =0. 
$

\subsection {Computations for higher dimensional hypersurfaces  } 
\noindent 

$\bullet$
\textbf {Singularities of Cubic Surface in $P^3$} 

The following list of singularities for cubic surfaces goes back to  Cayley. Schlafli was the first to classify the various types. We use the modern notation for the combination of singularities on a given surface, following Bruce and Wall \cite{BW}.

\textbf {$ A_1 : f= (x^2 + y^2 + z^2 + xy+ xz+ yz)w + 2xyz = 0$ } \\
$
S(t)=1+4t+6t^{2}+4t^{3}+(t^{4}+\ldots \\
\tau=1, \  \  ct=4, \  \  st=4, \  \  mdr=3, \  \  def=1, \  \  lc=4.$

\textbf {$ A_2 : f= f=(x+y+z)(x+2y+3z)w + xyz = 0$ } \\
$
S(t)=1+4t+6t^{2}+4t^{3}+2(t^{4}+\ldots \\
\tau=2, \  \  ct=3, \  \  st=4, \  \  mdr=2, \  \  def=1, \  \  lc=4.$

\textbf {$ 2A_1 : f=xzw+(z+w)y^2+x^3+x^2y+xy^2+y^3  = 0$ } \\
$
S(t)=1+4t+6t^{2}+4t^{3}+2(t^{4}+\ldots \\
\tau=2, \  \  ct=3, \  \  st=4, \  \  mdr=2, \  \  def=1, \  \  lc=4. $

\textbf {$ A_3 : f=xzw +(x+z)(y^2-x^2-z^2) = 0$ } \\
$
S(t)=1+4t+6t^{2}+4t^{3}+3(t^{4}+\ldots \\
\tau=3, \  \  ct=3, \  \  st=4, \  \  mdr=2, \  \  def=2, \  \  lc=3. $

\textbf {$ A_1+A_2 : f=x^3+y^3+x^2y+xy^2 + y^2z +xzw = 0$ } \\
$
S(t)=1+4t+6t^{2}+4t^{3}+3(t^{4}+\ldots \\
\tau=3, \  \  ct=3, \  \  st=4, \  \  mdr=2, \  \  def=2, \  \  lc=3.$

\textbf {$ A_4 : f= y^2z +yx^2-z^3+xzw = 0$ } \\
$
S(t)=1+4t+6t^{2}+4(t^{3}+\ldots \\
\tau=4, \  \  ct=3, \  \  st=3, \  \  mdr=2, \  \  def=3, \  \  lc=3.$

\textbf {$ 3A_1 : f=y^3+y^2(x+z+w)+4xzw; = 0$ } \\
$
S(t)=1+4t+6t^{2}+4t^{3}+3(t^{4}+\ldots \\
\tau=3, \  \  ct=3, \  \  st=4, \  \  mdr=2, \  \  def=2, \  \  lc=3. $

\textbf {$ 2A_2 : f= x^3+ y^3+x^2*y + xy^2 + xzw = 0$ } \\
$
S(t)=1+4t+6t^{2}+5t^{3}+4(t^{4}+\ldots \\
\tau=4, \  \  ct=2, \  \  st=4, \  \  mdr=1, \  \  def=1, \  \  lc=4.$

\textbf {$ A_1+A_3 : f=wxz+(x+z)(y^2-x^2) = 0$ }  \\
$
S(t)=1+4t+6t^{2}+4(t^{3}+\ldots \\
\tau=4, \  \  ct=3, \  \  st=3, \  \  mdr=2, \  \  def=3, \  \  lc=3.$

\textbf { $ A_5 : f=wxz + y^2z+x^3 - z^3  = 0$} \\
$
S(t)=1+4t+6t^{2}+5(t^{3}+\ldots \\
\tau=5, \  \  ct=2, \  \  st=3, \  \  mdr=1, \  \  def=1, \  \  lc=3.$

\textbf {$ D_4 : f=w(x+y+z)^2+xyz = 0$ } \\
$
S(t)=1+4t+6t^{2}+4(t^{3}+\ldots\\
\tau=4, \  \  ct=3, \  \  st=3, \  \  mdr=2, \  \  def=3, \  \  lc=3. $ 

\textbf { $ 2A_1+A_2 : f=wxz+y^2(x+y+z) = 0$ } \\
$
S(t)=1+4t+6t^{2}+4(t^{3}+\ldots\\
\tau=4, \  \  ct=3, \  \  st=3, \  \  mdr=2, \  \  def=3, \  \  lc=3. $

\textbf {$ A_1+A_4 : f=wxz+y^2z+yx^2 = 0$ } \\
$
S(t)=1+4t+6t^{2}+5(t^{3}+\ldots\\
\tau=5, \  \  ct=2, \  \  st=3, \  \  mdr=1, \  \  def=1, \  \  lc=3. $

\textbf { $ D_5 : f=wx^2+xz^2+y^2z;  = 0$}  \\
$
S(t)=1+4t+6t^{2}+5(t^{3}+\ldots\\
\tau=5, \  \  ct=2, \  \  st=3, \  \  mdr=1, \  \  def=1, \  \  lc=3. $

\textbf {$ 4A_1$  (Cayley surface) : $f=w(xy+xz+yz)+xyz= 0$ } \\
$
S(t)=1+4t+6t^{2}+4(t^{3}+\ldots \\
\tau=4, \  \  ct=3, \  \  st=3, \  \  mdr=2, \  \  def=3, \  \  lc=3. $ 

\textbf {$ A_1+2A_2 : f=wxz+xy^2+y^3  = 0$ } \\
$
S(t)=1+4t+6t^{2}+5(t^{3}+\ldots \\
\tau=5, \  \  ct=2, \  \  st=3, \  \  mdr=1, \  \  def=1, \  \  lc=3. $

\textbf {$ 2A_1+A_3 : f=wxz+(x+z)y^2  = 0$ }  \\
$
S(t)=1+4t+6t^{2}+5(t^{3}+\ldots \\
\tau=5, \  \  ct=2, \  \  st=3, \  \  mdr=1, \  \  def=1, \  \  lc=3. $

\textbf { $ A_1+A_5 : f=wxz +y^2z + x^3  = 0$}  \\
$
S(t)=1+4t+6(t^{2}+\ldots\\
\tau=6, \  \  ct=2, \  \  st=2, \  \  mdr=1, \  \  def=2, \  \ lc= \infty. $

\textbf {$ E_6 : f=wx^2 + xz^2 + y^3  = 0$ } \\
$
S(t)=1+4t+6(t^{2}+\ldots\\
\tau=6, \  \  ct=2, \  \  st=2, \  \  mdr=1, \  \  def=2, \  \ lc= \infty. $

\textbf {$ 3A_2 : f=wxz+ y^3  = 0$ } \\
$
S(t)=1+4t+6(t^{2}+\ldots\\
\Delta(t)=2t^{3}+5t^{4}+6(t^{5}+\ldots\\
\tau=6, \  \  ct=2, \  \  st=2, \  \  mdr=1, \  \  def=2, \  \ lc= \infty. $

$\bullet$
\textbf { Kummer quartic surface with 16 nodes } \\
$
f= x^4 + y^4 + z^4 - y^2z^2 - z^2x^2 - x^2y^2 - x^2w^2 - y^2w^2 - z^2w^2 + w^4 = 0  \\ 
S(t)= 1+4t+10t^2+16t^3+19t^4+16(t^5+\ldots \\
\tau= 16,  \  \ ct=5, \  \ st=5, \  \ mdr=3, \  \ def=6, \  \ lc=5.$

$\bullet$
\textbf { Octic nodal surface with 144 nodes } \\
$ f=16(x^8+y^8+z^8+w^8)+224(x^4y^4+x^4z^4+x^4w^4+y^4z^4+y^4w^4+z^4w^4)+
2688x^2y^2z^2w^2-9(x^2+y^2+z^2+w^2)^4 . \\
S(t)=1+4t+10t^2+20t^3+35t^4+56t^5+84t^6+116t^7+149t^8+180t^9+206t^{10}+224t^{11}+231t^{12}+224t^{13}+206t^{14}+
180t^{15}+158t^{16}+148t^{17}+145t^{18}+144(t^{19}+\ldots \\
\tau= 144, \  \  ct=15, \  \  st=19, \  \  mdr=9, \  \  def=9, \  \  lc=15. $

$\bullet$
\textbf { A non nodal octic surface (van Straten) } \\
$f=x^6y^2 - 2x^4y^4 + x^2y^6 - x^6z^2 - 3x^4y^2z^2 - 3x^2y^4z^2 - y^6z^2 + 5x^4z^4 + 10x^2y^2z^4 + 5y^4z^4 - 
8x^2z^6 - 8y^2z^6 - 508z^8 + 1024z^6w^2 - 640z^4w^4 + 128z^2w^6 - 8w^8.$ \\
$\tau(S/J)= 124, \tau(S/I)= 100$ so this surface is not nodal \\
S(t)=$1+4t+10t^{2}+20t^{3}+35t^{4}+56t^{5}+84t^{6}+116t^{7}+149t^{8}+180t^{9}+206t^{10}+224t^{11}+
231t^{12}+224t^{13}+206t^{14}+180t^{15}+157t^{16}+139t^{17}+128t^{18}+125t^{19}+124(t^{20}+\ldots $
invariants: $ct=15, \  \  st=20, \  \  mdr=9, \  \  def=8, \  \  lc=16. $ 

$\bullet$
\textbf {Singularities in $P^4$} 

\textbf{Cubic 3-fold with $10$ nodes }\\
$f=x_1^3+x_2^3 + x_3^3 + x_4^3+x_0x_1x_2 - x_0x_3x_4=0 $ \\
$
S(t)=1+5t+10(t^2+\ldots \\
\tau= 10,  \  \ ct=3, \  \ st=2, \  \ mdr=2, \  \ def=5,  \  \ lc= \infty .$

\textbf{Burkhardt quartic 3-fold with $45$ nodes }\\
$
f=x_0^4-x_0(x_1^3+x_2^3+x_3^3+x_4^3)+3x_1x_2x_3x_4 \\
S(t)=1+5t+15t^{2}+30t^{3}+45t^{4}+51t^{5}+45(t^{6}+\ldots \\
\tau=45, \  \  ct=6, \  \  st=6, \  \  mdr=4, \  \  def=15, \  \  lc=6. $

\textbf { Quintic with $125$ nodes} \\
$
f=x_0^5+x_1^5+x_2^5+x_3^5+x_4^5 -5x_0x_1x_2x_3x_4 = 0 \\
S(t)=1+5t+15t^{2}+35t^{3}+65t^{4}+101t^{5}+135t^{6}+155t^{7}+155t^{8}+135t^{9}+125(t^{10}+\ldots \\
\tau=125, \  \  ct=9, \  \  st=10, \  \  mdr=6, \  \  def=24, \  \  lc=9. $

$\bullet$
\textbf{Nodal cubic with $15$ nodes in $\PP^5$ } \\
$
f=4(x_0^3+x_1^3+x_2^3+x_3^3+x_4^3+x_5^3)-(x_0+x_1+x_2+x_3+x_4+x_5)^3=0 \\
S(t)=1+6t+15t^{2}+20t^{3}+15(t^{4}+\ldots \\
\tau =15, \  \  ct =4, \  \  st =4, \  \  mdr =3, \  \  def =9, \  \ lc=4. $

\section{Some properties for invariants} 

For any projective hypersurfaces, we compute Hilbert-Poincar\'{e} series and  some invariants, all the computations are made by Singular pachage.

If the hypersurfaces is nodal, Tjurina number $\tau $ is the number of nodes.

To decide if the hypersurface $f$ is nodal, we consider the following ideals: 
$J=J_f$ and $I=I_f$ the radical of $J$ and the invariants $\tau (S/J) $  and $\tau (S/I) $
The stability constant  $\tau (S/J) $ represent the total Tjurina number and for a nodal hypersurfaces represent the number of nodes.  
The stability constant  $\tau (S/I) $ represent the number of singularities and if we have $\tau (S/J) =\tau (S/I) $, the hypersurface $f=0$ are nodal.

Recall that Hilbert-Poincar\'{e} series of a graded $S$-module $M$ of finite type is defined by 
\begin{equation} 
\label{eq21B}
HP(M)(t)= \sum_{k\geq 0} (\dim M_k)t^k
\end{equation}
and that for smooth hypersurface we have 
\begin{equation} 
\label{eq31b}
F(t)=HP(M(f_s))=  \frac{(1-t^{d-1})^{n+1}}{(1-t)^{n+1}}=(1+t+t^2+\ldots + t^{d-2} )^{n+1}.
\end{equation}
In particular, if we set $T=T(n,d)=(n+1)(d-2)$, it follows that $M(f_s)_j=0$ for $j>T$ and $\dim M(f_s)_j=\dim M(f_s)_{T-j}$ for $0 \leq j \leq T$.
\medskip 

Smooth hypersurfaces have the same  Hilbert-Poincar\'{e} series associated to Fermat type

$F(t)=(1+t+t^2+\ldots + t^{d-2} )^{n+1} =\sum_{k=0}^{k=T}a_{k}t^{k}$.

In a very important recent paper \cite{H}, J. Huh has constructed for any homogeneous polynomial $h \in S_d$ a log-concave sequence $\mu^i(h)$, $i=0,...,n$ via the mixed multiplicities of the pair of ideals $J_h$ and $m=(x_0,...,x_n)$ in $S$ and showed surprising relation to the topology of the corresponding hypersurface $V(h):h=0$ in $\PP^n$.

It is natural therefore to ask whether our sequences $\dim M(f)_i$ give rise to some log-concave sequences. Here is the result.

\begin{prop} \label{prop1}
Let $f \in S_d$ and set as usual $T=(n+1)(d-2)$.

\begin{itemize}
\item If the hypersurface $V(f):f=0$ is smooth, then the sequence $a_k=\dim M(f)_k$, $k=0,...,T$ is log-concave.

\item If the hypersurface $V(f):f=0$ is singular, then the sequence $a_k=\dim M(f)_k$, $k=0,...,ct(V(f))$ is log-concave, but the extended sequence $a_k=\dim M(f)_k$, $k=0,...,ct(V(f)), ct(V(f))+1$ may be no longer log-concave even for nodal plane curves.
\end{itemize}

\end{prop}

\proof

For the smooth case, the coefficients $a_k$ appear as the coefficients of the polynomial
$$(1+t+t^2+...+t^{d-2})^{n+1}.$$
In other words, they are obtained by applying $n+1$ convolution products of sequences, starting with the constant sequence $1,1,...,1$. Since this constant sequence is clearly log-concave and the convolution product preserves the log-concavity, see \cite{H}, the first claim is proved.

The first part of the second claim is obvious by the definition of $ct(V(f))$. To justify the second part, consider the case where the curve is Lemniscate of Bernoulli . 

\endproof
 
\begin{rk}  \label{rem1}
If the hypersurface $V(h)$ has only isolated singularities, Huh's invariants $\mu^i(h)$ are very easy to compute: in fact $\mu^i(h)=(d-1)^i$ for $0 \leq i <n$ and 
$$\mu^n(h)= (d-1)^n - \mu(V(h)),$$
see Example 12 in \cite{H}. It follows that these invariants are not sensitive to the position of the singularities of $V(h)$ as our invariants $a_k=\dim M(h)_k$ are.
\end{rk} 

It is easy to see that one has 
\begin{equation} 
ct(D)=mdr(D)+d-2 \text { and }S(t)-F(t)=\defect(D)t^{ct(D)+1}+\ldots
\end{equation}
with
\begin{equation} 
 def(D)=\defect S_{(n+1)(d-2)-ct(D)-1}(\NN).
\end{equation} 

Note that computing the Hilbert-Poincar\'{e} series of the Milnor algebra $M(f)$ using an appropriate software like Singular, is much easier than computing the defects $\defect S_k(\NN)$, because the Jacobian ideal comes with a given set of $(n+1)$ generators $f_0,...,f_n$, while the ideal $I$ of polynomials vanishing on $\NN$ has not such a given generating set.

\section{Projectively rigid hypersurfaces} 

Recall that the automorphism group of the projective space $\PP^n$ is the linear group $G=PGL(n+1)$. It follows that any hypersurface $D$ in $\PP^n$ can be moved around by using translations by elements of $G$. A natural question is whether there are any other equisingular deformations of $D$ inside $\PP^n$. A recent answer to this question was given by E. Sernesi, who proved the following result, see formula $(5)$ and Corollary 2.2 in \cite{S} and refer to
\cite{S0} for the general theory of deformations.

\begin{prop} \label{sernesi}
Let $D:f=0$ be a degree $d$ hypersurface in $\PP^n$ having only isolated singularities. Let $\hat J$ be the saturation of the Jacobian ideal $J$ of $f$. Then 
the vector space $\hat J_d/J_d$ is naturally identified with the space of first order locally trivial deformations of $D$ in $\PP^n$ modulo those arising from the above $PGL(n+1)$-action.
\end{prop}

In view of this result we introduce the following definition.

\begin{definition}
\label{def2}
We say that a degree $d$ hypersurface $D:f=0$ is projectively rigid if $\hat J_d/J_d=0$. 
\end{definition}

\begin{ex}
\label{ex1}
Let $D$ be a Cayley cubic surface in $\PP^3$, i.e. a cubic surface having 4 $A_1$ singularities.
Since any two such surfaces differ by an element of the group $G=PGL(3)$, it follows that $D$ is
projectively rigid.
\end{ex}
In this section we show that many other classical surfaces and 3-folds are rigid in this sense. To do this, we first first explain how to compute the dimension of the vector space $\hat J_d/J_d$.
The formula (2.3) in \cite{DRR} gives
$$\dim S_d/\hat J_d=\tau(D)-\defect_d\Sigma_f.$$ 
Next we apply Theorem 1 in \cite{DRR} and get
$$\defect_d\Sigma_f=\dim M(f)_{T-d}-\dim M(f_{smooth})_d.$$
If we put everything together we get the following.

\begin{prop} \label{prop2}
$$\dim \hat J_d/J_d=\dim M(f)_{d}-\dim M(f_{smooth})_d+\dim M(f)_{T-d}-\tau(D).$$
\end{prop}

\begin{cor} \label{cor1}
A hypersurface $d$ degree $D:f=0$ in $\PP^n$ with isolated singularities and $d \leq ct(D)$ is projectively rigid
if and only if $ \dim M(f)_{T-d}=\tau(D)$.
\end{cor}

By a direct inspection of the results listed in section 2.2 we get the following.

\begin{cor} \label{cor2}
\noindent
\begin{itemize}

\item A cubic surface in $\PP^3$ with isolated singularities is projectively rigid if and only if it is one of the following types: 
$A_4$, $A_1+A_3$, $A_5$, $D_4$, $2A_1+A_2$, $A_1+A_4$, $D_5$, $4A_1$, $A_1+2A_2$, $2A_1+A_3$, $A_1+A_5$, $E_6$, $3A_2$ .

\item In $\PP^4$, the cubic type $10$ $A_1$, Burkhardt quartic type $45A_1$ and quintic type $125A_1$ are projectively rigid.

\item In $\PP^3$, the Kummer surface and octic type $144A_1$ are not projectively rigid.

\item Nodal cubic with $15$ nodes in $\PP^5$ is not projectively rigid.

\end{itemize}

\end{cor}

\section{Singular program to compute our invariants} 
For mathematical computations, we can use any CAS (Computer Algebra Systems) software like Mathematica, Matlab or Maple, but for Algebraic Geometry, the best are Singular, Macaulay2 or CoCoA.

Singular is a computer algebra system for polynomial computations, with special emphasis on commutative and non-commutative algebra, algebraic geometry, and singularity theory, developed at the University of Kaiserslautern (see http://www.singular.uni-kl.de/ and also \cite{GP}).\\
The following Singular program compute the Hilbert-Poincar\'{e} series and our invariants and decide if the hypersurface is nodal and projectively rigid.
\begin{verbatim} 
// run with command:  inv(n), where n+1 = number of variables 
proc inv(int n){
LIB "primdec.lib"; //library for  radical(Ideal);
ring R=0,(x(0..n)),dp;
poly f;
// you can change here the polynomial f
f=(x(0)^2+x(1)^2)^2-2*(x(0)^2-x(1)^2)*x(2)^2;  // Lemniscata 
int d=deg(f);
int T=(n+1)*(d-2); // T+1=maximal index of stabilization
ideal J=jacob(f);
ideal G=std(J);
ideal G2=std(radical(G));
int n1=mult(G), n2=mult(G2);
print("Total Tjurina= "+string(n1)); 
print("Number of singularities= "+string(n2)); 
if (n1==n2)
{ print("The hypersurface is nodal ! ");}
else
{print("The Hypersurface  is not nodal ! ");}
intvec v;
int i;
for (i=1; i<=T+2;i++) {v[i]=size(kbase(G,i-1));}
ring R=0,t,ds;
poly S; 
S=0;
for(i=1;i<=T+2;i++) {S=S+v[i]*t^(i-1); }
poly F=((1-t^(d-1))/(1-t))^(n+1); // Fermat series
poly D=S-F;
print("F(t)=Hilbert Poincare series for any non singular hypersurface ");
print(F);
print("S(t)=Hilbert Poincare series for singular hypersurface");
print(S);
number def1=leadcoef(D);
string stau1, sct1, sst1, smdr1, sdef1, slc1;
i=1;
while (coeffs(D,t)[i,1]==0) {i=i+1;}
int ct1=i-2;
int mdr1=ct1-(d-2);
i=T+2;
while (coeffs(S,t)[i,1]==coeffs(S,t)[i-1,1]) {i=i-1;}
int st1=i-1;
i=1;
while ((coeffs(S,t)[i,1]*coeffs(S,t)[i+2,1]<=coeffs(S,t)[i+1,1] *coeffs(S,t)
[i+1,1]) && (i<=1+st1))
{i=i+1;}
if (i==2+st1)
{ slc1=  "  lc="+"infinity ";}
else
{slc1=  "  lc="+string(i);}
stau1="   \\tau= "+ string(coeffs(S,t)[T+2,1]);
sct1 = "  ct="+string(ct1);
sst1 = "  st="+string(st1);
smdr1= "  mdr="+string(mdr1);
sdef1= "  def="+string(def1);
string ss1="invariants=  "+stau1+sct1+sst1+smdr1+sdef1+slc1;
print(ss1);
if ( coeffs(S,t)[d+1,1] - coeffs(F,t)[d+1,1] + coeffs(S,t)[T-d+1,1] 
- coeffs(S,t)[T+2,1] == 0)
{ print("The hypersurface is projectively rigid");}
else 
{ print("The hypersurface is not projectively rigid");}
};
\end{verbatim}


\end{document}